\documentclass[letterpaper, 10 pt, conference]{ieeeconf}  
\IEEEoverridecommandlockouts                              
\usepackage[math]{kurier}
\usepackage[sc]{mathpazo}

\usepackage{amssymb}
\usepackage{amsmath,amssymb}
\usepackage{amsbsy}
\usepackage{pgfgantt}
\usepackage{xcolor}
\usepackage{times}
\usepackage{soul}
\usepackage[utf8]{inputenc}
\usepackage[small]{caption}
\usepackage{microtype}
\usepackage{graphicx}
\graphicspath{ {figures/} }
\usepackage{wrapfig}
\usepackage{booktabs} 
\usepackage{subfiles}
\usepackage{csquotes}
\usepackage{cite} 
\usepackage{multicol}
\usepackage{latexsym} 
\usepackage{url}
\usepackage{lastpage}
\usepackage[sc]{mathpazo}
\usepackage{setspace}
\usepackage{graphicx}
\usepackage{enumerate}
\usepackage{float}
\usepackage[font=footnotesize]{caption}
\usepackage{pgfplots}
\usepgfplotslibrary{fillbetween}	
\usepackage{bbold}
\usepackage{comment}
\usepackage{mathtools, cuted}
\usepackage{authblk}
\usepackage{subcaption}
\usepackage{etaremune}
\usepackage{array}
\usepackage{tikz}
\usetikzlibrary{shapes}
\usetikzlibrary{calc} 

\usepackage[
  bookmarks=false,
  pdfpagelabels=false,
  hyperfootnotes=false,
  hyperindex=false,
  colorlinks,urlcolor= blue,citecolor=blue,
]{hyperref}

\usepackage{xcolor}
\definecolor{blue4}{HTML}{051A33}
\definecolor{blue3}{HTML}{153760}
\definecolor{blue2}{HTML}{275183}
\definecolor{blue1}{HTML}{3F689A}
\definecolor{blue}{HTML}{0029A3}
\definecolor{cfblue}{rgb}{0.258824,0.258824,0.435294}
\definecolor{dblue}{rgb}{.098,.243,.424}
\definecolor{dorange}{rgb}{0.72, 0.506, 0.125}
\definecolor{llblue}{rgb}{.447,.643,.831}

\newcommand{\balpha}{\boldsymbol{\alpha}}

\newcommand{\bbeta}{\boldsymbol{\beta}}

\newtheorem{theorem}{Theorem}

\newtheorem{remark}{Remark}

\newtheorem{definition}{Definition}
\newtheorem{proposition}{Proposition}

\DeclareMathOperator*{\argmin}{\arg\!\min} 
\DeclareMathOperator{\E}{\mathbb{E}}

\begin{document}
	
\title{\LARGE \bf
	Maximizing Road Capacity Using Cars that Influence People
}

\author[1]{Daniel A. Lazar}
\author[2]{Kabir Chandrasekher}
\author[1]{Ramtin Pedarsani}
\author[3]{Dorsa Sadigh}
\affil[1]{Electrical and Computer Engineering, University of California, Santa
    Barbara}
\affil[2]{Electrical Engineering, Stanford University}
\affil[3]{Computer Science \& Electrical Engineering, Stanford University}

	\maketitle
	\thispagestyle{empty}
	\pagestyle{empty}

	\begin{abstract}

The emerging technology enabling autonomy in vehicles has led to a variety of new problems in transportation networks, such as planning and perception for autonomous vehicles. Other works consider social objectives such as decreasing fuel consumption and travel time by platooning. However, these strategies are limited by the actions of the surrounding human drivers. In this paper, we consider \emph{proactively} achieving these social objectives by influencing human behavior through planned interactions. \emph{Our key insight is that we can use these social objectives to design local interactions that influence human behavior to achieve these goals.} To this end, we characterize the increase in road capacity afforded by platooning, as well as the vehicle configuration that maximizes road capacity. We present a novel algorithm that uses a low-level control framework to leverage local interactions to optimally rearrange vehicles. We showcase our algorithm using a simulated road shared between autonomous and human-driven vehicles, in which we illustrate the reordering in action.

\end{abstract}

\section{Introduction}
In recent years, the field of autonomous driving has experienced significant advances in the design of planning and perception algorithms~\cite{sadigh2016information,sadigh2016planning,Sadigh:EECS-2017-143,gray2013robust,raman2015reactive,vitus2013probabilistic,hermes2009long,vasudevan2012safe}. However, such techniques mainly consider a single autonomous vehicle without actually studying its societal effects such as its impact on commuter delay. On the other hand, recent work in transportation study how autonomous vehicles can achieve broader objectives such as increasing fuel efficiency and decreasing latency through techniques such as platooning~\cite{van2015fuel,adler2016optimal,lioris2017platoons,lazar2017routing,lazar2018routing,lazar2018altruistic}.
The mobility benefits of platooning autonomous or semi-autonomous vehicles have been largely studied in the literature for both freeway networks \cite{Shladover:1978lo,Yi:2006hb,Vander-Werf:2002fh,Van-Arem:2006ai,Arnaout:2014ul} and urban networks 	\cite{lioris2017platoons,Askari:2016fy}. Another line of work has focused on designing efficient scheduling policies at intersections by leveraging the autonomy of vehicles and their communication capabilities \cite{Lee:2012eu,Dallal:2013ph,Miculescu:2014pr,Colombo:2015df,Tallapragada:2015eb,Zhang:2016bc}. 
While previous techniques are quite valuable at solely studying local or global level properties, the two paradigms can actively influence each other.
This emphasizes the urgency and importance of understanding the interactions between humans and autonomous vehicles on shared roads as well as the global effects of these interactions.
For example, autonomous vehicles can stabilize traffic flow by damping shock waves in congested roads~\cite{wu2017emergent, wu2018stabilizing, cui2017stabilizing, stern2018dissipation}.
However, previous techniques do not \emph{proactively} create the circumstances that would enable them to positively impact these societal objectives: cars do not go out of their way to platoon, and they do not find congested traffic to stabilize.

\begin{quote}
 	\emph{Our key insight is that we need to actively influence human behavior by designing local interactions informed by desirable societal objectives for the shared road. }
 \end{quote}

\begin{figure}
 \centering
 \includegraphics[scale=0.25]{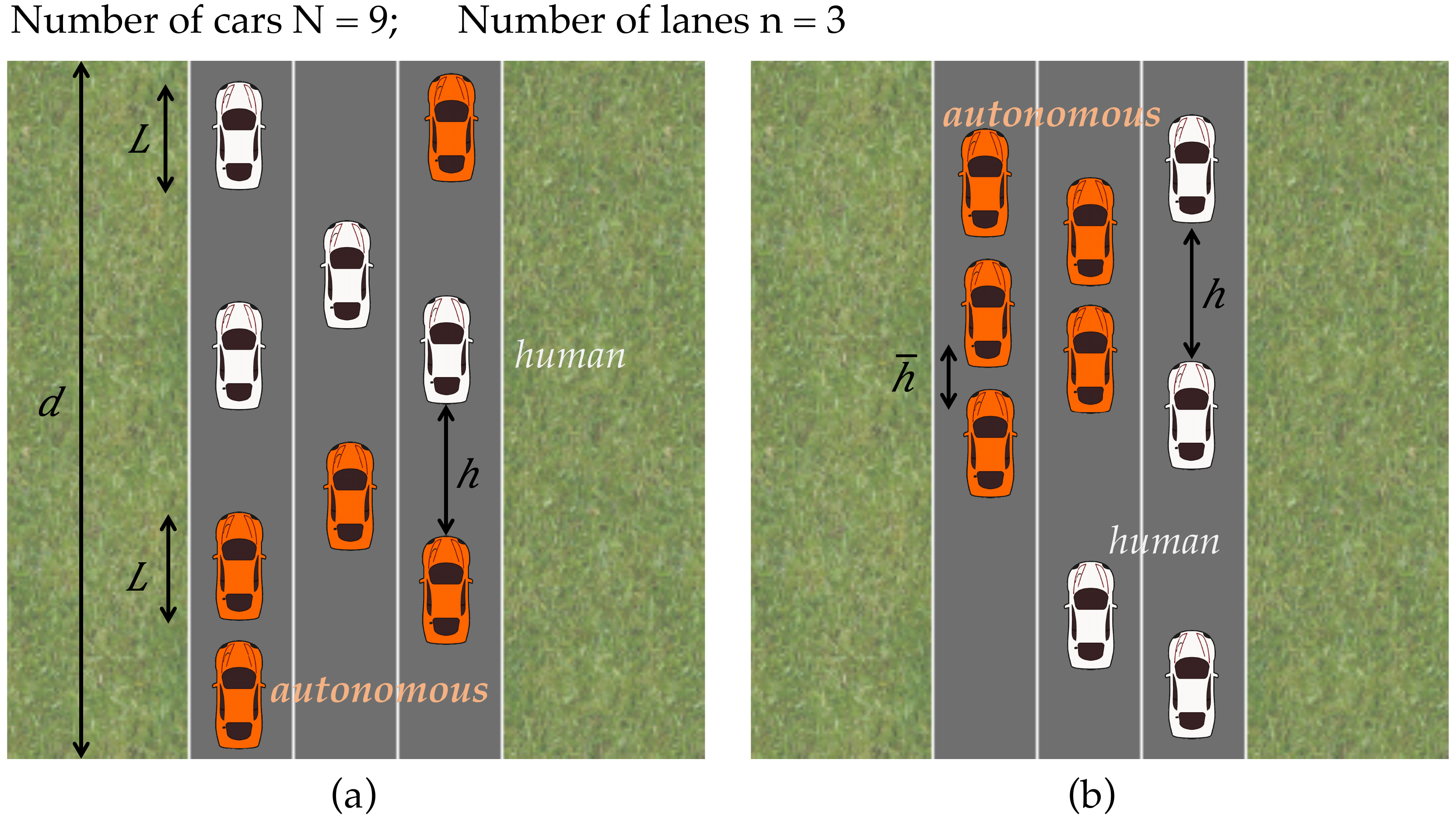}
 \caption{Road shared by autonomous (orange) and human-driven (white) cars.}	
 \label{fig:frontfig}
 \end{figure}

For instance, by anticipating a human's impatient response, an autonomous vehicle can influence a human driver to switch lanes by slowing down in front of them~\cite{sadigh2018planning}. 
Leveraging this insight, we plan to develop algorithms for \emph{interaction-aware} and \emph{socially-aware} autonomous vehicles that actively design their environments to reflect the most socially advantageous configurations.
Specifically, in this work, we design algorithms for autonomous cars to maximizes the road capacity. They achieve this via a sequence of interactions with  human drivers to reorder the vehicles on the road in order to platoon as much as possible. 
 For instance, as shown in Fig.~\ref{fig:frontfig}, if all autonomous cars (orange cars) on a road take interactive local actions that enforce human-driven cars (white cars) to open up some space, the autonomous cars can then leverage the opportunity to form a platoon (Fig.~\ref{fig:frontfig}~(b)). 
 
In this work, we build on notions of road capacity on shared roads~\cite{lazar2017routing} to include the idea of optimal lane allocation and vehicle ordering within a lane.  We then provide theoretical results on the structure and qualities of the optimal solution. We demonstrate that the optimal allocation and ordering can practically be achieved using state-of-the-art interaction-aware controllers~\cite{sadigh2016planning}. Our contributions in this paper are summarized as follows:
\begin{itemize}
	\item Leveraging local interactions between autonomous and human-driven cars to modify the road configuration in open loop to an optimal configuration.
	\item Precisely characterizing the optimal vehicle reordering and the price of no control, which quantifies the benefits of our algorithm.
	\item Implementing the controller for a mixed autonomy setting via simulation showing tight agreement with our theoretical results.
\end{itemize}
We defer proofs of all theoretical results to the appendix.
\section{Running Example}
In this section, we introduce a running example that describes an instance of the problem of our interest.
Imagine a road shown in Fig.~\ref{fig:frontfig} shared by autonomous and human-driven cars. The orange cars represent autonomous cars, and the white cars represent human-driven cars.

The vehicles, autonomous or human-driven, probabilistically arrive at any of the three lanes $1$, $2$, $3$ (where lane $1$ is the leftmost lane).
As we will discuss throughout the paper, the vehicles keep a headway between
each other, and this headway is much smaller if two autonomous cars are
following each other, i.e., $\bar{h} < h$. This reduced headway is described by the term \emph{platooning}, which intuitively means that autonomous cars,
unlike human-driven cars, are capable of coordinating with each other and driving
close to each other to save space and energy.

Our goal in this work is to start with a road configuration as shown in Fig.~\ref{fig:frontfig}~(a) and plan for the autonomous cars to navigate intelligently for the goal of converging to a configuration similar to that of Fig.~\ref{fig:frontfig}~(b).

We theoretically construct the most efficient road configuration based on an
optimal vehicle lane assignment and ordering problem and study its properties. In addition, to achieve this
optimal configuration, we leverage local interaction-aware controllers for
autonomous cars. As shown in previous work~\cite{sadigh2016planning}, the
autonomous vehicle can indirectly \emph{affect} a human-driven car to take
desirable actions that result in the ideal road configuration.  In constructing
these configurations, we assume that the heavy, local optimizations are solved
in a distributed fashion (separately by each vehicle), and that a
central controller exists which functions only to send and receive simple messages to coordinate actions among vehicles.

\section{Formalism}\label{sct:model}
Our goal in this paper is to leverage the low level interactions between autonomous and human-driven cars to arrange the vehicles in a fashion that is desirable from the road network's perspective.
To achieve this goal, we plan to answer two questions:
\begin{enumerate}
\item Assuming autonomous and human-driven vehicles can reorder themselves from any configuration to another (for example from Fig.~\ref{fig:frontfig}~(a) to Fig.~\ref{fig:frontfig}~(b)), what is the optimal configuration for the network?
\item How can we leverage the capability of autonomous cars to enforce this configuration?
\end{enumerate}

To study this problem, we first discuss road \emph{capacity} as a desirable property for the road network. Road capacity has been used as a common measure to be optimized for a desirable traffic network~\cite{Askari:2016fy,askari2017effect}. We model road capacity in mixed autonomy as in our previous works \cite{lazar2017routing,lazar2018routing}
	
	To model road capacity, we consider how many vehicles can be packed onto a road, with all vehicles traveling at nominal velocity and respecting their reaction time constraints. To this end, each car has a certain headway in front of it depending on the type of the car, as well as the type of the car it follows. In this model motivated by previous work in capacity modeling~\cite{Askari:2016fy, askari2017effect, lazar2017routing, lazar2018poa}, we assume that human-driven cars follow all vehicles at a distance of $h$ and autonomous cars follow other autonomous cars by a distance of $\bar{h}$ (typically $\bar{h} < h$) and follow human-driven cars with a distance of $h$ (as shown in Fig.~\ref{fig:frontfig}~(a)). Note that these quantities will vary by road as they depend on the nominal speed on that road. Further, the capacity on a road depends on the ordering of the vehicles on the road; we assume vehicle types are determined as the result of a Bernoulli process.
		
	We consider the capacity on a multi-lane road to be the sum of capacities of each lane. To find lane capacity, we first find the average headway, which is a function of the \emph{level of autonomy} of the lane, defined as $\alpha_i = \frac{y_i}{x_i + y_i} \in [0,1]$, where $x_i$ and $y_i$ are respectively the volume of human-driven and autonomous vehicles in lane $i$. For instance, the level of autonomy for the leftmost lane in Fig.~\ref{fig:frontfig}~(a) is $\alpha_1 = \frac{2}{4}$. 
	
	Let $h_j$ be the headway in front of vehicle $j$, and let $H_N$ be the total
    expected headway for $N$ vehicles on road $i$. Then, using linearity of expectation and the Bernoulli assumption, $H_N = \sum_{j=1}^{N-1}\E[h_j]=
    (N-1)(\alpha_i^2\bar{h} + (1-\alpha_i^2)h)$.
	
	If each vehicle has length $L$, then as the number of vehicles $N$ increases, the average space taken by a vehicle approaches $\alpha_i^2\bar{h} + (1-\alpha_i^2)h+L$. The capacity is the number of vehicles that can be fit on a lane of length $d$, so
	\begin{equation}\label{eq:cap2_phys_params}
	c(\alpha_i) = \frac{d}{\alpha_i^2\bar{h} + (1-\alpha_i^2)h+L} = \frac{d}{k_1 -\alpha_i^2 k_2} \; ,
	\end{equation}
	where we define two constants $k_1 \triangleq L + h$ and $k_2 \triangleq h - \bar{h}$, and $k_1>k_2$. Note that throughout the paper we consider same parameters, $d, h, \bar{h}, L$, for all the lanes.\\
	
	\noindent \textbf{Upper Bound on Capacity. }\label{sct:cap_UB}
	If we do not assume that the vehicles will arrive from some Bernoulli process and instead allow some arbitrary process to determine the ordering of the vehicles, the capacity model described above will no longer be accurate.  Though the specific form of the capacity function will depend on the ordering, we define a capacity model that serves as an upper bound on the capacity for any ordering.  This upper bound corresponds to the situation in which all autonomous cars are ordered optimally, meaning they are all adjacent and can form one long platoon.
	
	In this case, each autonomous car has a short headway $\bar{h}$ and each human-driven vehicle has headway $h$, yielding the average headway $\alpha_i \bar{h} + (1-\alpha_i)h$. This results in the following capacity function:
	\begin{align}\label{eq:capacity_UB}
	c^{UB}(\alpha_i) &= \frac{d}{\alpha_i \bar{h} + (1-\alpha_i)\bar{h} + L} = \frac{d}{k_1 -\alpha_i k_2} \; .
	\end{align}
	
	Since for all $i$ lanes $\alpha_i \in [0,1]$, it is clear that $c^{UB}(\alpha_i) \ge c(\alpha_i)$. Further, $c^{UB}$ serves as a tight upper bound with regards to vehicle ordering with a large number of vehicles.
	
\subsection{Characterizing Optimal Lane Assignment and Ordering}\label{sct:high_level_opt}
	We now pose optimization problems for vehicle lane choice when operating at capacity, under the capacity models in \eqref{eq:cap2_phys_params} and \eqref{eq:capacity_UB}. We also provide theoretical results pertaining to the optimal lane assignment and attendant total road capacities.
	
	\noindent \textbf{Optimal Lane Assignment and Ordering. } 
	\label{sct:routing_bernoulli}
	Assuming that the autonomy level for each lane $\alpha_i$ is chosen by a system operator, we now describe the high-level optimal lane assignment problem that maximizes the total capacity of the road.

Here we consider a multi-lane road to have a total capacity that is the sum of the capacity of each lane. 
All lane parameters will be the same, and the only difference between capacities in the lanes is due to different lane levels of autonomy $\alpha_i$.
	
	Let $\bar{\alpha}$ denote the overall autonomy level on the road, \emph{i.e.} the fraction all cars that are autonomous. Further, let $\balpha = \begin{bmatrix} \alpha_1 , \ldots , \alpha_n \end{bmatrix}^T$. Then, we define
	\begin{align}
	C(\balpha) &\triangleq \sum_{i=1}^n c(\alpha_i) \; , \label{eq:cap_shorthand} \\
	G(\balpha) &\triangleq \sum_{i=1}^n (\alpha_i - \bar{\alpha})c(\alpha_i) \; . \label{eq:constrnt_shorthand}
	\end{align}
	The social planner's optimization problem is as follows:  
	\begin{align}
	\balpha^* = &\arg\max_{\balpha} C(\balpha) \label{eq:maximization} \\
	&\text{s.t.} \; G(\balpha)=0 \; , \label{eq:constraint}
	\end{align}
	where $\alpha_i \in [0,1]$ for all lanes and \eqref{eq:constraint} constrains the solution to have overall autonomy level equal to $\bar{\alpha}$, the autonomy level of the traffic feeding the road. To see this, observe that $\sum_{i=1}^n \alpha_i c(\alpha_i)$ is the sum of autonomous cars on all lanes and has to be equal to $\bar{\alpha} \sum_{i=1}^n c(\alpha_i)$, which implies \eqref{eq:constraint}.
	Moreover, note that all the lanes have the same length, i.e. $d_i = d$;
    thus, we implicitly assume that lane utilization is at capacity. 

	\begin{theorem}\label{thm:single_mixed_lane}
		Consider an optimization problem of the form \eqref{eq:maximization}. Any solution will have at most one lane with mixed autonomy, \emph{i.e.} at most one lane with $\alpha_i \in (0,1)$.
	\end{theorem}
	\vspace{2mm}
	\begin{remark}
		This theorem formalizes the intuition that mixing autonomous and regular cars only decreases capacity, as it lessens the likelihood of adjacent autonomous vehicles, the condition necessary for platooning.
	\end{remark}
	
	Now that we have information that narrows down the set of possible solutions, we can derive a closed-form expression for optimal lane assignment, as expressed in the following theorem. 
	\begin{theorem}\label{prop:optimal_routing}
		Let $\balpha^* = \begin{bmatrix} \alpha^*_1, \; \ldots, \;  \alpha^*_n \end{bmatrix}^T$ be an optimum of \eqref{eq:maximization}, with autonomy levels ordered in decreasing order. Let $m$ denote the last lane with full autonomy, \emph{i.e.} $m=\max i$ s.t. $\alpha^*_i=1$, where $m=0$ implies that there are no lanes with full autonomy. Then, $$m = \lfloor  \frac{\bar{\alpha}n(k_1-k_2)}{k_1-\bar{\alpha}k_2} \rfloor.$$
	\end{theorem}
    \vspace{2mm}
	
	\begin{remark}
		One can also derive a closed-form expression for $\alpha^*_{m+1}$ by solving a quadratic equation, but for the sake of brevity we omit this result.
	\end{remark}
		After characterizing the optimal solution to \eqref{eq:maximization}, we now describe the lane assignment that minimizes total capacity. This will be used to characterize the cost of declining to optimally route vehicles, which is developed in Section \ref{sct:cost_of_no_planning}.
	\begin{proposition}\label{prop:min_capacity}
		Let $\balpha_*$ denote the worst-case lane assignment, corresponding to the minimum sum of capacities:
		\begin{equation*}
		\balpha_* = \argmin_{\balpha}C(\balpha) \; \text{s.t. \eqref{eq:constraint}.}
		\end{equation*}

		Then,
			$\balpha_* = \mathbb{1}_n^T \bar{\alpha}$ and 
			$C(\balpha_*) = nc(\bar{\alpha})$,
		where $\mathbb{1}_n$ is the row vector of ones of length $n$.
	\end{proposition}
	\vspace{2mm}
	See Appendix \ref{pf:min_capacity} for a sketch of the proof. Intuitively, the worst-case lane assignment is when all lanes have autonomy level equal to the overall autonomy level. This means that any perturbation from uniform autonomy level allows benefits from platooning.
	
	\noindent \textbf{Upper Bound on Capacity Under Optimal Ordering. } \label{sct:routing_UB}
		Using the upper bound capacity in \eqref{eq:capacity_UB}, we now formulate an optimization problem to find the theoretical maximum capacity on a road with $n$ lanes, with autonomy level $\bar{\alpha}$. 
	
	For reasons of readability, we use $\bbeta$ when discussing autonomy levels that are solutions to the upper bound capacity maximization. As before, let $C^{UB}(\bbeta) \triangleq \sum_{i=1}^n c^{UB}(\beta_i)$ and $G^{UB}(\bbeta) \triangleq \sum_{i=1}^n (\beta_i - \bar{\alpha})c^{UB}(\beta_i)$. The optimization problem is then:
	\begin{align}
	&\max_{\bbeta} C^{UB}(\bbeta) \label{eq:maximization_UB} \\
	&\text{s.t.} \quad G^{UB}(\bbeta) = 0 \; , \label{eq:constraint_UB} \\
	& \qquad \beta_i \in [0,1] \; \forall \; i \nonumber
	\end{align}        
	
	We next present a lemma showing that the total capacity when using the upper
    bound on capacity is invariant to lane assignment.
	\begin{proposition}\label{prop:UB_equal_cost}
		Any feasible solution to \eqref{eq:maximization_UB} has cost $nc^{UB}(\bar{\alpha})$.
	\end{proposition}
	
	This proposition follows from the ordering-invariant nature of the capacity upper bound. One can also construct a proof similar to that of Proposition \ref{prop:min_capacity}.
	
	\begin{remark}\label{rem:UB_routing}
		Given that all feasible lane assignments have the same total capacity, there are two notable lane assignments to consider. One is the uniform assignment, $\bbeta^*=\begin{bmatrix} \bar{\alpha}, \; \ldots, \; \bar{\alpha} \end{bmatrix}^T$. Another is is the assignment with a maximum of one mixed lane, such as the one developed in Theorem \ref{prop:optimal_routing} for the original capacity function. In this case, there will be same number of purely autonomous lanes, which is apparent from the proof of the theorem. However, the level of autonomy in the mixed lane will differ from that in the solution to \eqref{eq:maximization}. Again, for brevity we do not state the closed-form solution for this autonomy level, but it can be found by manipulating the constraint \eqref{eq:constraint_UB}.
	\end{remark}

	\subsection{Local Interaction on Shared Roads}
	\label{sct:low_level_opt}
Now that we have discussed optimal vehicle arrangement, we address
the second question of our interest, i.e., \emph{how can we leverage the
capability of autonomous cars to enact this optimal arrangement?}

The solution of the optimal lane assignment problem provides an optimal level of autonomy $\balpha^*$ based on the model of road capacity as we have discussed so far. 
Since the cars are driving on shared roads, we leverage the power of autonomous cars to enforce this level of autonomy by allowing them to navigate intelligently on the road.

The autonomous cars in each lane can take actions that \emph{affect} the behavior of the human-driven cars. For instance, they can cause humans to change lanes, slow down, or speed up.
These local actions can in fact create a reordering of the vehicles.
Our goal is to intelligently create such reorderings in order to get closer to the optimal capacity for the road. In this section, we discuss how to initiate this local reordering.

We model the local interactions between one autonomous car and one human-driven car as a dynamical system: $x^{t+1} = f(x^t, u^t_A, u^t_H)$, where the states of the world $x^t$ evolve based on the actions of the human-driven car $u_H^t$ and the autonomous car $u_A^t$.
Our goal in these local interactions is to design a controller for the autonomous car, i.e., $u_A^*(u_H^*)$, that not only achieves reaching its destination but also can affect the actions of the human-driven car. In other words, it can target specific reactions from the human-driven car.
We follow the work of Sadigh et al.~\cite{sadigh2016planning} to design such controllers by formulating the problem as a nested optimization:
\begin{equation}
u_A^* = \arg \max_{u_A} R_A(x, u_A, u_H^*(x, u_A)) \; .	
\label{eq:optimal_robot}
\end{equation}	
Here $R_A$ denotes the autonomous car's reward function which, if optimized, outputs a sequence of actions for the autonomous car $u_A^*$.
We assume this reward function is a combination of objectives such as avoiding collisions, keeping distance to the road boundaries, and the vehicle reaching its goal. Specific goals, such as influencing humans to switch lanes, also appear in $R_A$.

The optimization in~\ref{eq:optimal_robot} depends on a model of the human driving behavior that outputs $u_H^*$.
In this work, we model humans as agents who optimize their own reward function:
\begin{equation}
u_H^* = \arg \max_{u_H} R_H(x,u_A^*, u_H) \; .
\end{equation}
Similarly, this reward function encodes the goals and objectives of the human-driven car such as collision avoidance or keeping a certain velocity and heading. The specific formulation of both $R_H$ and $R_A$ can be found in~\cite{Sadigh:EECS-2017-143}.
The reward $R_H$ is usually learned from demonstration via techniques such as inverse reinforcement learning (IRL)~\cite{abbeel2004apprenticeship,levine2012continuous,ziebart2008maximum}, in which trajectories of human driving are collected in an offline setting and $R_H$ is then computed based on the collected trajectories.

We note that optimizing the autonomous car's reward function $R_A$ (as in \eqref{eq:optimal_robot}) can implicitly \emph{manipulate} the actions of the human-driven car. Similarly, the actions of the human-driven car can influence the autonomous vehicle's actions. To avoid infinite regress from this endless recursion, we approximate this interaction using a Stackelberg (leader-follower) game~\cite{osborne1994course}, where we assume the autonomous car plays first influencing the human driver, while the human driver \emph{observes} actions taken by the autonomous vehicle. In practice, this approximation accurately reflects observed behavior since the autonomous vehicles replans repeatedly using a model predictive controller~\cite{sadigh2016planning}.

Though in theory this controller can be extended to autonomous vehicles influencing multiple humans in the presence of other autonomous vehicles, the nested optimization in \eqref{eq:optimal_robot} is very computationally intensive. We therefore limit our scheme to have each autonomous vehicle take into account its influence over a single human driver, while treating the actions of other drivers, as well as the actions of other autonomous vehicles, as fixed with respect to the autonomous vehicle's actions. 

Through this interaction between an autonomous car and a human-driven car~\eqref{eq:optimal_robot}, we can intelligently design reward functions for the autonomous car $R_A$, which then leads to actions from the autonomous car $u_A^*$ that help the vehicle navigate for more efficient road usage by affecting the actions of the human-driven car $u_H^*$. We refer to the human influenced as being \emph{paired} with the autonomous vehicle.

	\subsection{Cost of lack of planning}\label{sct:cost_of_no_planning}
We have now discussed how one can solve the optimal lane assignment problem, as well as how local interactions affect the actions of the human-driven car.
Before presenting our simulation results, we first present theoretical limits on how much performance can be improved by optimally assigning lanes and optimally reordering in the mixed lane. 
This will allow us to gauge the efficacy of our control scheme, as well as judge when it is worthwhile to attempt more difficult vehicle manipulation. We do this by examining the magnitude of potential gains that can be achieved by these more complicated maneuvers.
To this end, we introduce the following two quantities:
	\begin{definition}\label{def:price_of_negligence}
		The \emph{price of negligence}, denoted $\Lambda$, is the maximum ratio between road capacity at optimal vehicle arrangement on a road to capacity at the worst-case lane assignment. This is formalized as follows:
				\begin{align*}
			\Lambda &= \max_{\bar{\alpha}, n, L, h,\bar{h}} \frac{C^{UB}(\bbeta^*)}{C(\balpha_*)} \\
			& \qquad \; \text{s.t. feasible network parameters.} 
		\end{align*}

	\end{definition}
	\begin{remark}
		Note that worst-case lane assignment in the denominator above maintains the Bernoulli assumption, \emph{i.e.} not considering worst-case car ordering. Under worst-case ordering, autonomous and human-driven cars would be interleaved, unlike the assumption leading to the capacity model in \eqref{eq:cap2_phys_params}. 
	\end{remark}

	\begin{definition}\label{def:price_of_no_control}
		The \emph{price of no control}, denoted $\Gamma$, is the maximum ratio between road capacity at optimal vehicle ordering on a road and road capacity at best-case vehicle lane assignment without reordering. This is formalized as
				\begin{align*}
			\Gamma &= \max_{\bar{\alpha}, n, L,h,\bar{h}} \frac{C^{UB}(\bbeta^*)}{C(\balpha^*)} \\
			& \qquad \; \text{s.t. feasible network parameters.}
		\end{align*}

	\end{definition}

	\noindent With this in mind, we present our bounds on these quantities.
	\begin{theorem}\label{thm:price_negligence_no_control}
		The \emph{price of negligence} is bounded by 
		\begin{align*}
			\Lambda \le 2\frac{L + h- \sqrt{(L + h)(L+\bar{h})}}{h-\bar{h}} \le 2 \; ,
		\end{align*}
		and the \emph{price of no control} is bounded by
		\begin{align*}
			\Gamma \le \frac{2n(L + h)}{(2n-1)(L + h) + \sqrt{(L + h)(L + \bar{h})}} \le \frac{2n}{2n-1} \; .
		\end{align*}
	\end{theorem}

	\begin{remark}
		The derived upper bounds on $\Lambda$ and $\Gamma$ are achieved by setting vehicle length $L$ and short headway $\bar{h}$ to 0. Intuitively, the gain due to platooning increases as the space taken up by a platoon decreases.
	\end{remark}
	\begin{remark}
		Note that $\Gamma \le \Lambda$, with equality when $n=1$. This is because when there is only one lane, there is no gain from optimal lane assignment.
	\end{remark}
	
	The implications of the upper bounds on price of negligence and price of no control are perhaps surprising. This means that if vehicle type is decided by a Bernoulli process, optimally ordering the vehicles will result in at most a factor of 2 increase in capacity. In a more realistic scenario, consider $L = 4$, $h = 30$, $\bar{h} = 11$. Then, $\Lambda \le 1.202$, meaning that the maximum possible improvement is approximately $20\%$.
	
	Note that a $20\%$ increase in capacity may yield far more than a $20\%$ improvement in road latency. For example, a fourth-order polynomial is commonly used in traffic literature to describe the relationship between vehicle flow on a road and the road's latency. Consider, as in \cite[Ch.~13]{sheffi1985urban}, road latency of the form $\ell(x) = t^0[1 + \rho(\frac{ x}{c'})^\sigma]$, where $x$ is the volume of traffic on the road, $c'$ is the ``practical capacity", $t^0$ is the free-flow travel time, and $\rho$ and $\sigma$ are typically $0.15$ and $4$, respectively. Using this model, under a given traffic flow, a $20\%$ increase in capacity leads to a $50\%$ reduction in latency due to congestion. In high congestion, this results in a large decrease in total latency.

	\section{Achieving Optimal Configuration}\label{sct:mid_level_opt}
	
	In this section we describe our algorithm for approaching the capacity increases yielded from the optimal lane assignment and optimal vehicle ordering. We describe our algorithm and discuss the theoretical upper bound for performance at each stage. Below is a summary of the policy:
	
	\begin{itemize}
		\item Phase 0: Cars are initialized in an intermixed configuration and ordering in all lanes is determined by a Bernoulli process with parameter $\bar{\alpha}$.
		\item Phase 1: Autonomous vehicles follow optimal lane assignment.
		\item Phase 2: Autonomous vehicles \emph{influence} human drivers in local interactions to follow optimal lane assignment.
		\item Phase 3: Autonomous vehicles start platooning in the mixed lane, if one exists.
	\end{itemize}

	Note that if we were to terminate the policy at phase 2, we would want the lane assignment to follow the solution to \eqref{eq:maximization}. However, if we terminate the policy at phase 3, the optimal lane assignment is the solution to \eqref{eq:maximization_UB}. The number of purely autonomous lanes in the solutions are the same (see Remark \ref{rem:UB_routing}), but the autonomy level in the mixed lanes will differ. We choose the optimistic assignment, which assumes that we can optimally rearrange the vehicles.

	\subsection{Phase 1: Autonomous cars follow lane assignment}\label{sct:phase1}
	
	\begin{figure}
		\centering
		\includegraphics[width=\columnwidth]{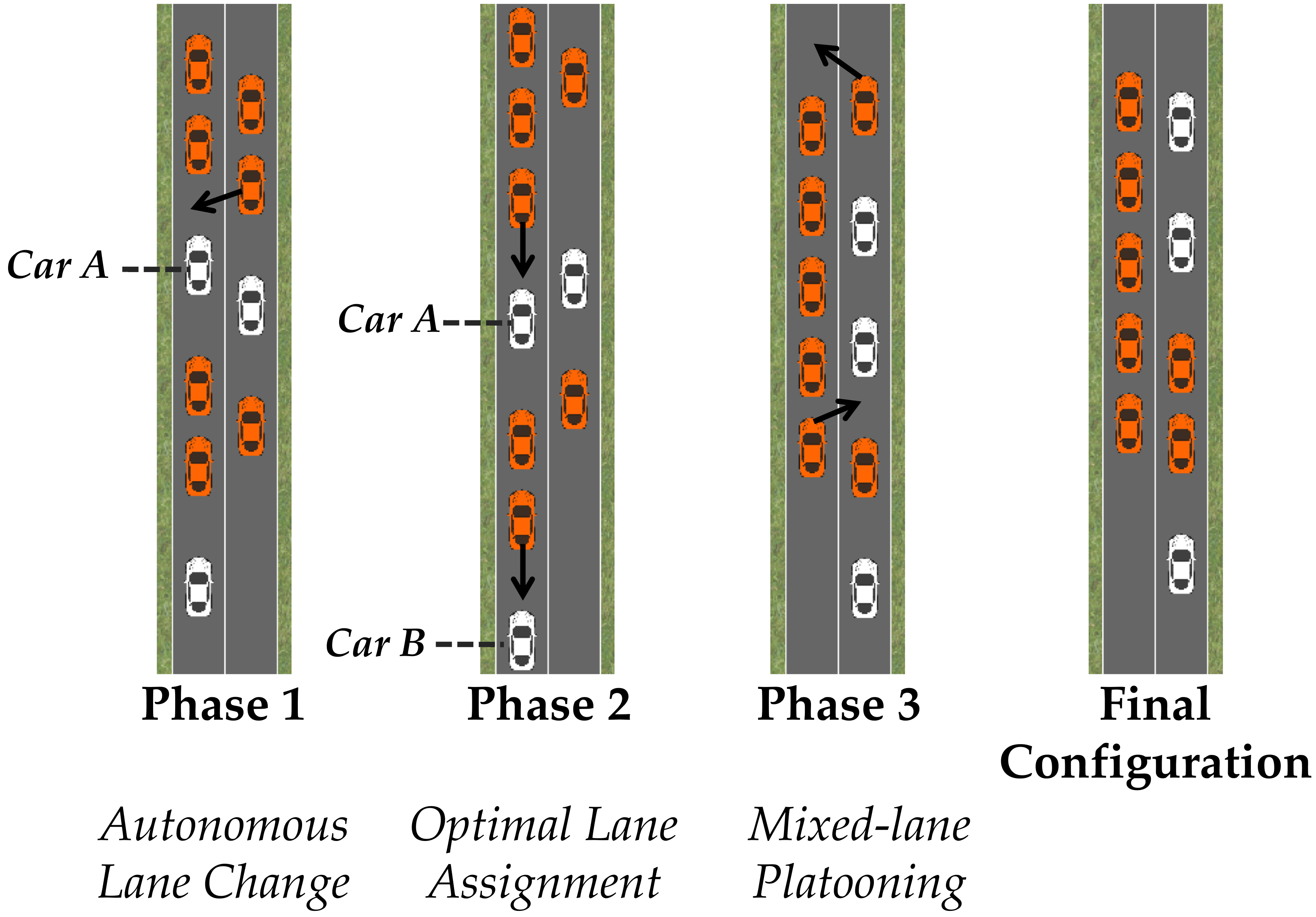}
		\caption{\footnotesize{Phase 1: Autonomous vehicles follow optimal lane assignment. Here, the acting autonomous car pairs with human-driven car A. Phase 2: Autonomous vehicles influence humans to follow optimal lane assignment. Here, the acting autonomous cars pair with the human-driven cars A \& B. Phase 3: Autonomous vehicles platoon in the mixed lane.}~~}
		\label{fig:phases}
	\end{figure}
	
	In this stage, pictured in Fig. \ref{fig:phases} (Phase 1), autonomous vehicles switch lanes to satisfy their optimal lane assignment determined by $\balpha^*$, the ordered solution to \eqref{eq:maximization_UB}. The human-driven vehicles are not assumed to follow the optimal lane assignment by the end of this stage.
	
	To accomplish this, each autonomous vehicle can be directed to switch to a specific lane, if there is a strong centralized control with close sensing and communication with the vehicles.  This can be accomplished in a decentralized manner as well, if the autonomous vehicles are given a vector of probabilities $\boldsymbol{q}^* = [q_1^* \; q_2^* \; \text{...} \; q_n^*]^T$ corresponding to $\boldsymbol{\alpha}^*$, where $q_i^*$ is the probability that an autonomous vehicle is assigned to lane $i$. The vehicles then sample from a generalized Bernoulli distribution with parameter $\boldsymbol{q}^*$ to determine which lane to inhabit. 
	
	However the desired lane is selected, in the process of switching lanes, each autonomous vehicle will determine its action by performing the nested optimization in~\eqref{eq:optimal_robot} while including behavior of the human-driven vehicles in their path. This pairing allows the autonomous car to model the human's responses, allowing it to merge even when there is not a large space for it, by relying on trailing human-driven vehicles slowing down in response to the autonomous car.
	
	\subsection{Phase 2: Autonomous cars influence human lane choice}\label{sct:phase2}
	
	This phase takes place only if the solution to \eqref{eq:maximization} yields at least one lane with all autonomous vehicles, \emph{i.e.} $\alpha_1 = 1$. Otherwise, all influence over human-driven cars occurs in the mixed lane, which takes place during Phase 3.
	
	If there are some lanes that under lane assignment have only autonomous vehicles, in this stage the autonomous vehicles in those lanes influence the human-driven vehicles to enforce this solution, as shown in Fig.~\ref{fig:phases} (Phase 2). Lane-by-lane, starting with lane 1, autonomous cars influence human drivers to leave their lanes. This is done by having each autonomous car check if there is a human-driven vehicle behind it; if so, it pairs with that vehicle and with the goal to influence them to merge rightward. In a decentralized manner, each autonomous vehicle in an autonomous lane monitors its surroundings, and if there are no human-driven vehicles to the left of it, it checks for a human-driven vehicle behind it. If there is one, it influences it to move to the lane on the right.
	
	By the end of this phase, lanes designated as purely autonomous will be so, with all autonomous vehicles in these lane forming long platoons. We expect the total capacity at the end of this phase to compare to the solution of \eqref{eq:maximization}.
	
	\subsection{Phase 3: Autonomous cars platoon in the mixed lane}\label{sct:phase3}
	
	Though the theoretical maximum capacity under the Bernoulli assumption has
    been achieved at this point, if there is a lane with mixed autonomy, the
    capacity can be further improved by platooning the vehicles in that lane. If
    the mixed lane is the only lane with autonomous vehicles, i.e., $\forall i,
    \; \alpha_i^*<1$, then autonomous vehicles can pair with human-driven
    cars behind them, leading them out of the lane, until they reach another
    autonomous car and platoon with it.
	
	However, if there is a lane with full autonomy, this scheme can be more easily achieved through a series of swaps, as seen in Fig.~\ref{fig:phases} (Phase 3). If there is a lone autonomous vehicle, it can merge into the adjacent full-autonomy lane, and another car can exit the lane to join a platoon in the mixed lane. The exiting vehicle can do this by pairing with the human-driven vehicle behind the vehicle it wants to platoon with, so that it will not be too cautious in switching lanes. Alternatively, the vehicle to be platooned with can slow down, and the merging vehicle will move in front of it. 
	
	By the end of this phase, total capacity will be lower bounded by the optimal capacity in \eqref{eq:maximization}, and upper bounded by the optimum of \eqref{eq:maximization_UB}.

     \begin{figure}[t]
 \centering

 \includegraphics[width=0.9\columnwidth]{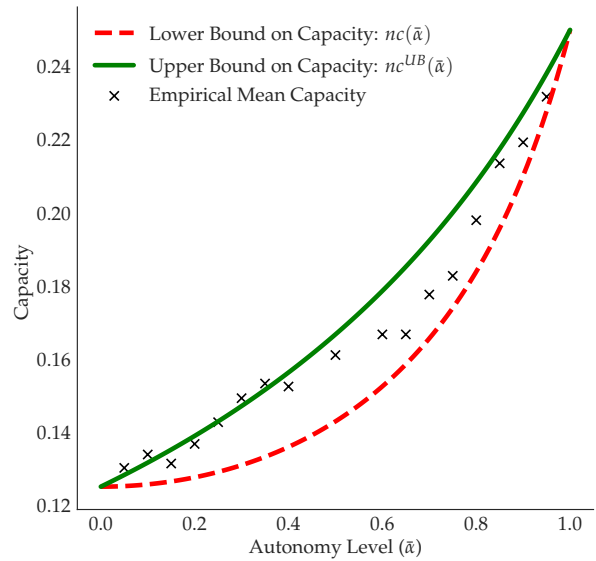}
 \caption{\textbf{Capacity with Local Interactions}. The achieved capacity at
 various autonomy levels plotted against the achievable capacity with optimal
 local interactions (green) and baseline capacity without any local
         interaction.}
 \label{fig:capacity_fig}
 \end{figure}

      \begin{figure}[t]
 \centering
 \includegraphics[scale=0.45]{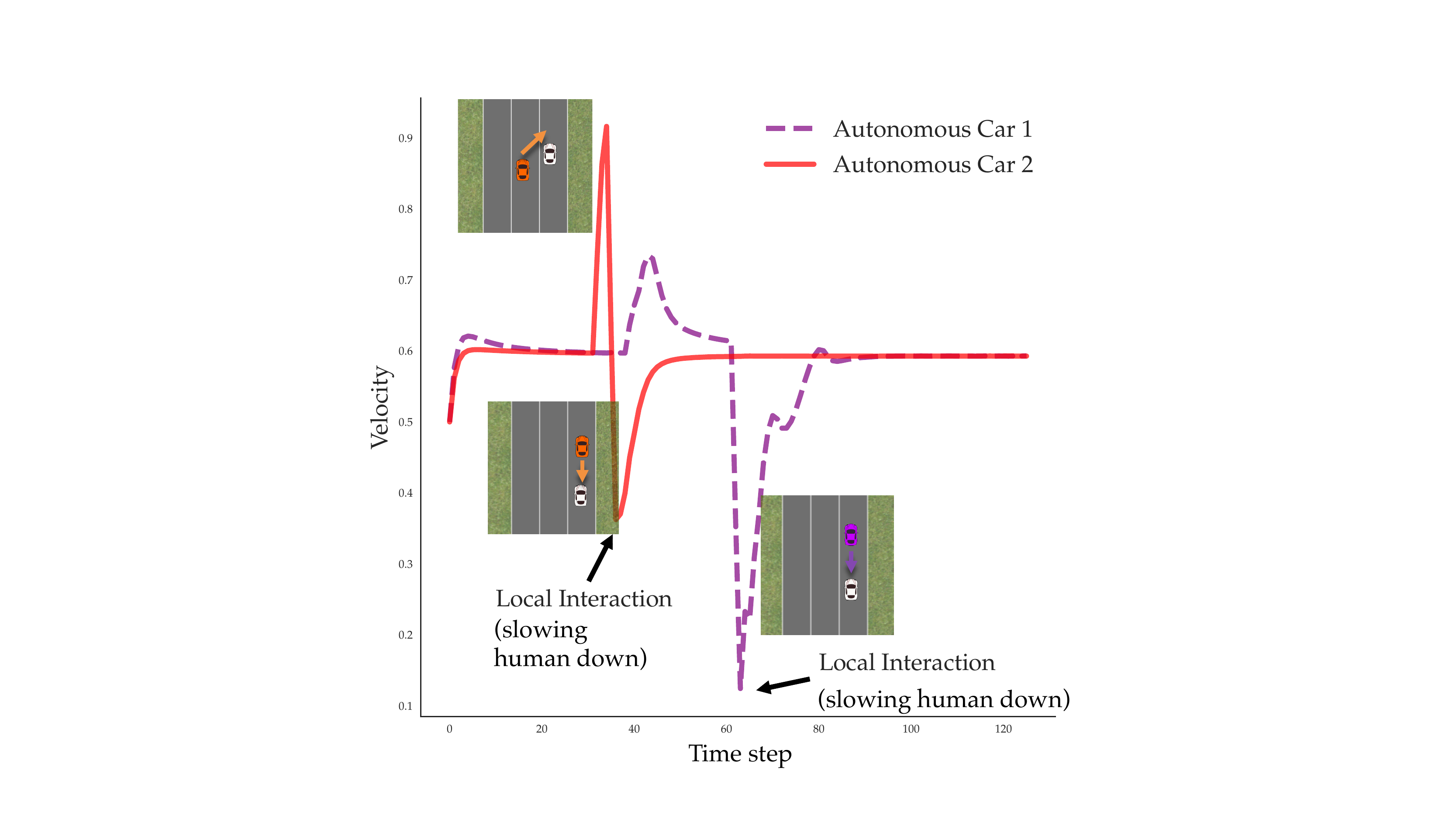}
 \caption{\textbf{Timing of Local Interactions}. The variation in velocity of two autonomous vehicles (purple and red), each influencing a human to merge lanes. The spikes indicate an autonomous vehicle accelerating to move in front of a human, and the dips indicate the autonomous vehicles slowing down to influence the human to merge away.}
 \label{fig:velocity_fig}
 \end{figure}

	In this section, we discuss our autonomous driving simulation framework and compare the results of the simulation to the theoretical bounds established in Section \ref{sct:high_level_opt}.
	
	\subsection{Experimental Setup}
	
We simulate a two-lane road with 20 vehicles, at varying levels of autonomy. To
determine vehicle arrangement at the start of the simulation, we create a
population with desired number of vehicles of each type, shuffle it, then
randomly place the vehicles into lanes, ensuring that each lane has the same
number of vehicles. This begins phase 0, in which cars are intermixed as the
result of a Bernoulli process.  Human-driven cars are simulated as agents
optimizing a reward function, specifically the rewards
in \cite{sadigh2016planning}.

	\subsection{Driving Simulator} 
	We use a simple point-mass model for the dynamics of each vehicle, where the state of the system is: $x = [x \enskip y \enskip \theta \enskip v]^\top$. $x,y$ represent the coordinates of the vehicle on the road, $\theta$ is its heading, and $v$ is its velocity. We assume two control inputs $u = [u_1 \enskip u_2]$ each representing the steering angle and the acceleration of the vehicle. We then represent the dynamics of each vehicle as the following, where $\gamma$ is the friction coefficient:
\begin{equation*}
[\dot{x} \enskip \dot{y} \enskip \dot{\theta} \enskip \dot{v}] = [v\cdot\cos(\theta) \quad v\cdot\sin(\theta) \quad v\cdot u_1 \quad u_2-\gamma\cdot v].
\end{equation*}

\subsection{Simulation Results}
Recall the situation of Fig.~\ref{fig:frontfig} and note that our goal is to
start in a configuration similar to the left side and end in a situation such as
the right.  Noting our theoretical upper and lower bounds on the effects of
lane choice and ordering, we show in Fig.~\ref{fig:capacity_fig} the results of the mid-level optimization
and compare it to the earlier derived upper and lower bounds. Constrained the computationally intensive nested optimization \eqref{eq:optimal_robot}, we simulate $200$ time steps using a setup with $20$ cars at varying autonomy levels. As expected,
each run lies between the two curves.  We note that the few data points which
exceed the upper bound are because the upper bound includes the headway of the frontmost vehicles -- this is a good approximation for a large number of vehicles, but the simulation includes only $20$. Note also that at lower
levels of autonomy, the data points are closer to the upper bound, whereas at
higher levels, the data points are closer to the lower bound.  This is due to
the fact that the nested optimization is done one car at a time and with only
one other car in the network. We illustrate the timing of the local
interactions in Fig.~\ref{fig:velocity_fig}, noting that when the velocity of
the autonomous cars drops dramatically, it indicates a movement intended to
influence a human-driven car. 

\subsection{Implementation Details}
In our implementation of the nested optimization, we used the software package Theano~\cite{bergstraal:2010-scipy,Bastien-Theano-2012} to symbolically compute all Jacobians and Hessians. Theano optimizes the computation graph into efficient C code, which is crucial for real-time applications. 
In our implementation, each step of our optimization is solved with a horizon
length $N=5$.  We run the large-scale simulations with $20$ cars in the network on
a cluster using $4$ CPU's with a maximum utilization of $36$ GB RAM between
them.  We make the code publicly available here:
\href{https://github.com/Stanford-HRI/MultilaneInteractions}{https://github.com/Stanford-HRI/MultilaneInteractions}.

	\section{Discussion}\label{sec:conclusions}
	\noindent \textbf {Summary. } We introduce a procedure to connect
    societal goals for shared roads to low-level interactions designed for autonomous vehicles to influence human drivers. To this end, we provide theoretical guarantees on the potential of optimally assigning lanes and re-ordering to maximize road capacity via platooning. Finally, we implement an algorithm to achieve this reordering.
	
	\noindent \textbf {Limitations. }
    The nested optimization is extremely computationally heavy, and is natural
    to be run in a distributed setting in which each car runs its own nested
    optimization.  It would be interesting to implement a fully distributed
    version of the algorithm and compare the results to the current
    implementation.  Additionally, we emphasize that the simulations were not
    run with real humans, but rather with agents whose reward functions were
    learned from humans. Each human is, however, different and this could have large impacts on the performance of the algorithm. We plan to address these shortcomings in future work.

	\noindent \textbf {Conclusion. } We demonstrate the connection between
    societal objectives and vehicle-level control in transportation networks and design a unifying scheme utilizing both to achieve optimal efficiency in the network.

\section*{Acknowledgement}\label{sec:acknowledgements}
Some of the computing for this project was performed on the Sherlock
cluster. We would like to thank Stanford University and the Stanford
Research Computing Center for providing computational resources and support that contributed to these research results. This work was supported in part by NSF grant no. CCF-1755808 and the UC Office of the President, grant no. LFR-18-548175.

	\section{Appendix}
	
	\subsection{Proof of Theorem \ref{thm:single_mixed_lane}}\label{pf:single_mixed_lane}
	Assume, for the purpose of contradiction, that there exists a solution $\boldsymbol{\alpha}^*$ to \eqref{eq:maximization} with components $\alpha_i^*, \; \alpha_j^* \in (0,1)$, where $i \neq j$.
	Assume without loss of generality that $\alpha_i \le \alpha_j$. If $\alpha_i = \alpha_j$, we can construct a new solution with same total capacity in which $\alpha_i \neq \alpha_j$. This is because $\frac{\partial }{\partial \alpha_i}C(\balpha) =\frac{\partial }{\partial \alpha_j}C(\balpha)$ and $\frac{\partial }{\partial \alpha_i}G(\balpha) =\frac{\partial}{\partial \alpha_j}G(\balpha)$, so the autonomy level on road $i$ can be decreased by some infinitesimal value $\epsilon$ and autonomy level on road $j$ increased by $\epsilon$, thereby maintaining the same capacity sum and satisfying the constraint.
	
	Now that we have established an optimal solution with $\alpha_i^* < \alpha_j^*$ with $\alpha_i^*, \; \alpha_j^* \in (0,1)$, we explore a further perturbation of our solution. We construct a new routing $\hat{\balpha}$ with $\hat{\alpha}_i = \alpha_i^* - \epsilon_1$ and $\hat{\alpha}_j = \alpha_j^* + \epsilon_2$, where $\epsilon_1, \; \epsilon_2 > 0$ with these unequal infinitesimal perturbations designed such that the constraint remains satisfied. Note that this can be accomplished with nonnegative perturbations as $\frac{\partial}{\partial \alpha_i}G(\balpha) = \frac{k_2 \alpha_i^2 - 2k_2 \alpha_i \bar{\alpha} + k_1}{(k_1 -k_2\alpha_i^2)^2}>0$, which is positive since $k_1>k_2$.
	
	Since we choose the relative scaling of $\epsilon_1$ and $\epsilon_2$ based on their contribution to the constraint function, we can find the relative change in the objective due to the perturbations by the expression	
	\begin{align*}
		&\frac{\frac{\partial}{\partial \alpha_j} C(\boldsymbol{\alpha})}{\frac{\partial}{\partial \alpha_j} G(\boldsymbol{\alpha})} - \frac{\frac{\partial}{\partial \alpha_i} C(\boldsymbol{\alpha})}{\frac{\partial}{\partial \alpha_i} G(\boldsymbol{\alpha})} \\
		&=\frac{2k_2 \alpha_j}{k_2 \alpha_j^2 - 2k_2\alpha_j \bar{\alpha} + k_1} - \frac{2k_2 \alpha_i}{k_2 \alpha_i^2 - 2k_2\alpha_i \bar{\alpha} + k_1} >0
	\end{align*}
	since $\frac{2k_2 \alpha_i}{k_2 \alpha_i^2 - 2k_2 \alpha_i \bar{\alpha} + k_1} > 0$ and
	\begin{align*}
		\frac{\partial}{\partial \alpha_i} \frac{2k_2 \alpha_i}{k_2 \alpha_i^2 - 2k_2\alpha_i \bar{\alpha} + k_1} &= \frac{2k_2(k_1-\alpha^2k_2)}{(k_2 \alpha_i^2 - 2k_2\alpha_i \bar{\alpha} + k_1)^2} > 0 \; ,
	\end{align*}
	as $k_1>k_2>0$. Therefore, the addition of $\epsilon_2$ to $\alpha_j$ increases the objective more than the subtraction of $\epsilon_1$ from $\alpha_i$ decreases it. Then, $C(\boldsymbol{\hat{\alpha}})>C(\balpha^*)$, so $\balpha^*$ is not an optimum, proving the theorem. \hfill \QED 
	
	\subsection{Proof of Theorem \ref{prop:optimal_routing}}\label{pf:optimal_routing}	
	Assume that $\bar{\alpha}<1$, so $m<n$. Let $g(\alpha_i) = (\alpha_i - \bar{\alpha})c(\alpha_i)$. Then using \eqref{eq:constraint},
	\begin{align*}
		0 &= G(\balpha^*) = \sum_{i=1}^{n}g(\alpha_i^*) = \sum_{i=1}^{m}g(1) + g(\alpha_{m+1}^*) + \sum_{i=m+2}^{n}g(0) \\
		&= mg(1)+g(\alpha^*_{m+1}) + (n-m-1)g(0) \; .
	\end{align*}
	We can solve for $m$, giving
	\begin{align*}
		m=\frac{(k_1-k_2)(n\bar{\alpha} (k_1-\alpha^*_{m+1})k_2 - \alpha^*_{m+1}(k_1 - \alpha^*_{m+1}\bar{\alpha} k_2))}{(k_1-\alpha^{*2}_{m+1})(k_1-\bar{\alpha} k_2)} \; .
	\end{align*}
	
	One can show that $\frac{\partial m}{\partial \alpha_{m+1}^*}<0$. Note that by definition, $\alpha^*_{m+1} \in [0,1)$ (as $m$ represents the index of the first lane such that $\alpha_i \neq 1$), so we can bound $m$ using these values:
	\begin{align*}
	\frac{\bar{\alpha} n(k_1-k_2)}{k_1-\bar{\alpha} k_2}-1 < m \le \frac{\bar{\alpha} n(k_1-k_2)}{k_1-\bar{\alpha} k_2} \; .
	\end{align*}
	
	Therefore, $m = \lfloor  \frac{\bar{\alpha} n(k_1-k_2)}{k_1-\bar{\alpha} k_2} \rfloor$. Note that in the case excluded at the start of the proof, when $\bar{\alpha} = 1$, this expression yields $m=n$, which is correct. Therefore, this expression is true for $\bar{\alpha} \in [0,1]$. \hfill \QED 
	
	\subsection{Proof of Proposition \ref{prop:min_capacity}}\label{pf:min_capacity}
	
	First note that given feasible routing vector $\balpha$, if there exists one element $\alpha_i>\bar{\alpha}$, then there must also exist element $\alpha_j<\bar{\alpha}$. This follows from the fact that $G(\mathbb{1}_n^T \bar{\alpha})=0$ and $\frac{\partial}{\partial \alpha_i}G(\balpha)>0$, as shown in the proof of Theorem~\ref{thm:single_mixed_lane}.
		
	Now we can prove the proposition by recursion. The base case is $\balpha = \mathbb{1}_n^T \bar{\alpha}$. For any other feasible routing, there are at most $n$ lane autonomy levels $\alpha_i$ not equal to $\bar{\alpha}$. Pick any pair of lanes $i$ and $j$ such that $\alpha_i>\bar{\alpha}>\alpha_j$ (which, when not in the base case, are guaranteed to exist by the fact above). Then, using the reverse of the mechanism in the proof of Theorem \ref{thm:single_mixed_lane}, we keep the constraint satisfied while decreasing $\alpha_i$ and increasing $\alpha_j$, while monotonically decreasing the road capacity. This continues until either $\alpha_i=\bar{\alpha}$, $\alpha_j=\bar{\alpha}$, or both. Now we have a routing vector with lower road capacity than the original, with at most $n-1$ lane autonomy levels not equal to $\bar{\alpha}$, proving the proposition.
	
	\subsection{Proof of Theorem \ref{thm:price_negligence_no_control}}\label{pf:price_negligence_no_control}
	To prove the first statement, note that due to Proposition \ref{prop:UB_equal_cost}, the cost for any feasible routing for the upper bound capacity function is equivalent to $nc^{UB}(\bar{\alpha})$. For given network parameters $k_1$ and $k_2h$, we find our price of negligence:
	\begin{align*}
		\Lambda \le  \max_{\bar{\alpha} \in [0,1] }\frac{nc^{UB}(\bar{\alpha})}{nc(\bar{\alpha})} = \max_{\bar{\alpha} \in [0,1] } \frac{k_1-\bar{\alpha}^2 k_2}{k_1-\bar{\alpha} k_2} \; .
	\end{align*}
	
	This term is concave with respect to $\bar{\alpha}$ for $\bar{\alpha} \in [0,1]$, with maximum at $\bar{\alpha} = \frac{k_1-\sqrt{k_1^2-k_1k_2}}{k_2}$. Plugging this in,
	\begin{align}
		\Lambda 
		&\le 2\frac{k_1-\sqrt{k_1^2-k_1k_2}}{k_2} = 2\frac{L + h- \sqrt{(L + h)(L+\bar{h})}}{h-\bar{h}} \nonumber \\ \label{eq:thm3-proof}
		&\leq 2.
	\end{align}
	To observe \eqref{eq:thm3-proof}, note that $k_1 - 2k_2 \leq \sqrt{(k_1 - k_2)k_1}$.
	
	For the second statement, first note that if $\bar{\alpha}=1$, the price of no control would be 1. Since $\Gamma \ge 1$, let us exclude that case. Now, assuming $\bar{\alpha} < 1$ we apply Theorem \ref{thm:single_mixed_lane} and Proposition \ref{prop:UB_equal_cost}:
	\begin{align}
		\Gamma & \le \frac{nc^{UB}(\bar{\alpha})}{\sum_{i=1}^{m}c(1) + c(\alpha^*_{m+1})+ \sum_{i=m+2}^{n}c(0)} \nonumber \\
		&\le \frac{nc^{UB}(\bar{\alpha})}{c(\alpha^*_1) + (n-1)c(0)} \label{eq:pfthm2ineq1} \; ,
	\end{align}
	where \eqref{eq:pfthm2ineq1} follows from the fact that $\Gamma \ge 1$ and $c(1) > c(\alpha^*_{m+1}) \ge c(0)$ so the denominator is minimized when $m=0$ (meaning there are no purely autonomous lanes).

	We can then solve $G(\balpha^*) = 0$ to find an expression for $\bar{\alpha}$ as a function of $\alpha^*_1$, yielding
	\begin{align*}
		\Gamma \le \frac{n(k_1 - \alpha^*_1 k_2)}{n(k_1-\alpha^*_1 k_2) - \alpha^*_1 k_2(1-\alpha^*_1)}
	\end{align*}
	We then consider bounding $\Gamma^{-1}$, which we show to be convex.
	\begin{align*}
		\frac{\partial^2}{\partial \alpha^{*2}_1}(1/\Gamma) &= \frac{2 k_2 (k_1^2 + 3k_1k_2(-1+\alpha^*_1)\alpha - k_2^2\alpha^{*3}_1)}{n(k_1-\alpha^{*2}_1k_2)^3}
	\end{align*}
	The inner term on the numerator is convex for $\alpha^*_1 \in [0,1]$, with second derivative $6k_2(k_1-\alpha^*_1 k_2)$. Its minimum over that interval occurs at $\alpha^*_1 = \frac{k_1-\sqrt{k_1^2-k_1 k_2}}{k_2}$, yielding $\frac{\partial^2}{\partial \alpha^{*2}_1}(1/\Gamma)>0$. Solving $\frac{\partial}{\partial \alpha^*_1}(1/\Gamma)=0$, we find that the minimum of this outer equation also occurs at $\alpha^*_1 = \frac{k_1-\sqrt{k_1^2-k_1 k_2}}{k_2}$. Using this,
	\begin{align*}
		\Gamma &\le \frac{4k_1n^2 - 2n(k_1+\sqrt{k_1^2 - k_1 k_2})}{4k_1n(n-1)+k_2} \\
		&= \frac{2n(L + h)}{(2n-1)(L + h) + \sqrt{(L + \bar{h})(L + \bar{h})}} \\
		&\leq \frac{2n}{2n-1}\; .
	\end{align*}
	As $k_2$ approaches $k_1$, the above inequality becomes tight. \hfill \QED

\end{document}